\documentclass[11pt]{amsart}
\usepackage{amsmath,amsthm,amssymb,amsfonts}
\usepackage{graphicx,color,epsfig,subfigure}
\usepackage{url,verbatim,multirow}

\title{A generalization of the methods of Brass, Harborth, and Nienborg}
\author{Brendon Stanton}
\address{Department of Mathematics, Iowa State University, Ames, Iowa 50011}

\begin{document}
\begin{abstract}
In 1995, Brass, Harborth and Nienborg disproved a conjecture of Erd\H{o}s when they showed that a $C_4$-free subgraph of the hypercube, $Q_n$,  can have at least $(\frac 12 +\omega(1))e(Q_n)$ edges.  In this paper, we generalize the idea of Brass, Harborth and Nienborg to provide good constructions of $Q_3$-free subgraphs of $Q_n$ for some small values of $n$.
\end{abstract}
\maketitle

\section{Main Idea}
The idea is to generalize the idea of Brass, Harboth and Nienborg~\cite{Brass1995} and apply it to other graphs.  We start with a base graph $G_k\subset Q_k$ which is $H$-free as well as a colored $Q_m$ graph.  We then split $G_k$ into $G_k^{(1)}$ and $G_k^{(2)}$ with a (not necessarily complete) parallel set of edges in between them.  We then consider our colored graph $Q_m$.  Since it is bipartite, we replace each vertex in one partite set with $G_k^{(1)}$ and the other with $G_k^{(2)}$.  We then use the colors of the edges in $Q_m$ to somehow determine which edges to include in our new $G_{k+m-1}\subset Q_{k+m-1}$ which will be $H$-free, assuming our algorithm works correctly.

\section{The BHN construction}

We start with an $aeo$-coloring of $Q_m$.  All $a$-colored edges are replaced with a copy of $P$.  For the $e$- and $o$-colored edges, we note that $G_k^{(i)}$ is again bipartite. We divide each of these into bipartite sets and add all edges in one bipartite set wherever we see an $e$-colored edge and all edges in the other set when we see an $0$-colored edge.  The new graph will be $C_4$-free.

It is interesting to note that assuming $G_k$ is maximal (but not necessarily maximum), the resulting construction will also be maximal (I think).

\subsection{Why it works}

We wish to try and find a $C_4$ in our newly constructed graph. Clearly, such a $C_4$ cannot be a subgraph of $G_k^{(i)}$ since $G_k^{(i)}$ is $C_4$ free.

Next, suppose our $C_4$ contains one edge in some $G_k^{(1)}$ and another in some $G_k^{(2)}$ with two edges going in between them.  If the two edges came from an $a$-labeled edge in $G_m$, then notice that $G_k^{(1)}\cup G_k^{(2)}\cup P$ forms a copy of $G_k$ and so there is no $C_4$.  If it comes from an $e$- or $o$-colored edge, then we see that since the two edges are adjacent in $G_k^{(i)}$, they must come from different partite sets and so one of the edges in our $e$- or $o$- edge is not present in our new graph.

Finally, the only case left is that all 4 of our vertices come from different copies  of $G_k^{(i)}$.  Then our each of our vertices are in the same partite set and this graph forms a 4-cycle in $Q_m$ so it contains at least one $e$-edge and one $o$-edge.  But since these contain edges from different partite sets, one of them is absent.

\section{First Idea}

We wish to copy the idea of Brass, Harborth and Nienborg exactly.  The only difference is that we now require two things from our $aeo$-coloring of $G_m$:
\begin{enumerate}
  \item Each $Q_3$ contains at least one $e$-edge \emph{and} one $o$-edge.
  \item Each $Q_2$ contains at least one $e$-edge \emph{or} one $o$-edge.
\end{enumerate}

To see that this creates a $Q_3$-free subgraph isn't much harder than before.  \\~\\

(1) No $Q_3$ can be contained entirely inside a $G_k^{(i)}$.\\~\\

(2) If our $Q_3$ is contained within a single edge of $Q_m$
then if that edge was an $a$-edge, we once again have a copy of $G_k$.  If it's an $e$- or $o$-edge, it is actually missing not one by two of the edges in between.\\~\\

(3) If our $Q_3$ is contained within a $C_4$ of $Q_m$, then it has at least one $e$- or $o$- edge.  Since there are two possible edges running between those in our new graph and they come from different partite sets, one of them is not in our new graph.\\~\\

(4) Finally, the only case is that our $Q_3$ is part of a $Q_3$ in $G_m$, but as in the last case, this contains one $e$- and one $o$-edge and so one of these is not an edge.\\~\\

\subsection{Why this doesn't work}

Since ex$(Q_3,Q_2)$=9, each cube in our graph must contain at least 3 $e$- or $o$-edges.  Thus, our $G_m$ must be at least 1/4 non-$a$ edges.  Thus, at every step we are omitting at least 1/8 of the possible edges within these edges.  This means that we would need to find a parallel class with at least 5/6 of the possible edges to stay above 3/4 of the total edges in the graph.

\subsection{Salvaging Something?}
Although this doesn't work asymptotically, we may be able to use it for small values of $n$.  Furthermore, since the construction is not likely to be maximal (see case 2 from before) we may be able to go back and add edges by hand.  This should be a much simpler job for my java program.

We can $aeo$-color a $Q_3$ by removing the edges $[00*]$, $[*11]$ and $[1*0]$ and assigning two of them color $e$ and one of them color $o$.  Let $G_k$ be a $Q_3$-free subgraph of $Q_k$.  Let $e_k$ be the number of edges in this graph and $p_k$ be the number of edges in a parallel class. (We wish to choose the parallel class with the largest number of edges or we get nothing).

Using our technique, we can create a $Q_3$ free $G_{k+2}$.  When creating this, we will have exactly 4 $G_k^{(1)}$s and 4 $G_k^{(2)}$s.  The total number of edges in these is $4(e_k-p_k)$.  We will be adding $9p_k$ edges corresponding to our $a$-colored edges and $3\cdot 2^{k-2}$ edges corresponding to our $e$- and $o$-colored edges.  This gives the recurrence $$e_{k+2} = 4e_k+5p_k + 3\cdot 2^{k-2}.$$

For $k=5$, we may use the unique construction given by Offner\cite{Offner2011} which has 72 out of 80 possible edges and a parallel class with all 16 edges.  Plugging this into the equation above gives $e_7=392$.  Since there are 448 total edges in $Q_7$, using that notation this gives $c(Q_3,7)\le 56$, which improves at least on the 1993 bound by Graham, Harary, Livingston and Stout\cite{Graham1993} of 62.

Applying this to $Q_4$ gives us $c(Q_3,6)\le 24$, which is worse than the exact result of 22.  In a $Q_3$-free $Q_6$ omitting 22 edges, we must have one parallel class with $\lfloor 22/6\rfloor =3$ omitted edges.  So if we take this construction we have $e_6 = 192-22=70$ and $p_6=32-3=29$ which gives $e_8=873$ and so $$c(Q_3,8)\le 1024-873=151.$$  I have no idea how this compares to known bounds.

\subsubsection{A simpler construction}
Actually, an easier way to do this is to just take a $C_4$ with one edge labeled $e$ and the other 3 labeled $a$.  This gives the recurrence $$e_{k+1} = 2(e_k-p_k)+3p_k+2^{k-2} = 2e^k + p_k + 2^{k-2}.$$

Plugging in our result for $e_6$ above with $p_6=29$, this gives $e_7=395$ and so $c(Q_3,7)\le 53$.  Then $p_7\ge 57$ and so $c(Q_3,8)\le 145$.

Edit:  This actually gives $e_7=385$, so the result is worse.

\subsubsection{A more elaborate construction}
Another $aeo$-coloring of $Q_4$ is possible using the following construction:

$$\begin{array}{|c|c|c|c|}
  \hline
   \text{$e$ edges} & \begin{array}{cc}
                \text{[*000]} \\
                \text{[*111]} \\
                \text{[1*01]} \\
                \text{[0*10]}
              \end{array}
   &
  \text{$o$ edges} &  \begin{array}{cc}
                \text{[00*1]} \\
                \text{[11*0]} \\
                \text{[101*]} \\
                \text{[010*]}
              \end{array} \\
  \hline
\end{array}$$ and the rest of the edges are colored with $a$.  This yields the recursion: $$e_{k+3} = 8e_k + 16p_k + 2^{k+1}.$$

Also note at this point that we may replace $e_k$ with the number of non-edges in $G_k$ and $p_k$ with the number of non-edges in our parallel class and the recursive formula remains the same.  Recall that $e_4=3$ and $p_4=0$ (using the number of omitted edges) so this gives $e_7=8\cdot 3 + 32 = 56$ gives us the same upper bound for $c(Q_3,7)$ as before.  However, for $k=5$ we have $e_k=8$ and $p_k=0$ so we get $c(Q_3,8)\le 128$ which is a substantial improvement on the bound above.  This also gives $c(Q_3,9)\le 352$.

I ran the $G_7$ construction through my program which verified that it was indeed $Q_3$-free with omitted edges:
\begin{verbatim}
[*000000],[000*010],[0000*11],[*000011],[000*101],[*000101],
[*000110],[0001*00],[*001001],[*001010],[*001100],[*001111],
[001000*],[0*10001],[0*10010],[00101*0],[0*10100],[0*10111],
[0*11000],[00110*1],[0*11011],[0*11101],[001111*],[0*11110],
[010000*],[01001*0],[01010*1],[010111*],[011*010],[0110*11],
[011*101],[0111*00],[100000*],[10001*0],[10010*1],[100111*],
[101*010],[1010*11],[101*101],[1011*00],[11*0001],[110*010],
[11*0010],[1100*11],[11*0100],[110*101],[11*0111],[1101*00],
[11*1000],[11*1011],[11*1101],[11*1110],[111000*],[11101*0],
[11110*1],[111111*]
\end{verbatim}

I tried to run a perturbation algorithm on this graph to see if we could perhaps do slightly better by first removing 2 edges and then adding edges to the resulting graph, but after 20 hours, my computer spat out the same graph.  To try first deleting 3 edges, I estimate that it would take my machine a bit over 2 years.

We may also give the more general recurrence.  Since there are $k$ parallel classes and $e_k$ total edges in our $G_k$, we may find $p_k\le \lfloor e_k/k\rfloor$ and so $$e_{k+3} \le 8e_k + 16\lfloor e_k/k\rfloor + 2^{k+1}.$$

Here is a table of bounds up to the point where we start getting a proportion of edges larger than 1/4, taking the lower bound from~\cite{Offner2011}:

$$\begin{array}{|c|cc|cc|}
\hline
&&&&\\
k & \text{LB} & \text{UB} & \text{LB}/e(Q_k) & \text{UB}/e(Q_k) \\
&&&&\\
\hline
&&&&\\
7 & 52 & 56 & \frac{52}{448}\approx 0.11607 & \frac{56}{448}= 0.125 \\
&&&&\\
8 & 119 & 128 & \frac{119}{1,024}\approx 0.11621 & \frac{128}{1,024}= 0.125 \\
&&&&\\
9 & 268 & 352 & \frac{268}{2,304}\approx 0.11632 & \frac{352}{2,304}\approx 0.15278 \\
&&&&\\
10 & 596 & 832 & \frac{596}{5,120}\approx 0.11641 & \frac{832}{5,120}\approx 0.16250 \\
&&&&\\
11 & 1,312 & 1,792 & \frac{1,312}{11,264}\approx 0.11648 & \frac{1,792}{11,264}\approx 0.15909 \\
&&&&\\
12 & 2,863 & 4,464 & \frac{2,863}{24,576}\approx 0.11650 & \frac{4,464}{24,576}\approx 0.18164 \\
&&&&\\
13 & 6,204 & 10,032 & \frac{6,204}{53,248}\approx 0.11651 & \frac{10,032}{53,248}\approx 0.18840 \\
&&&&\\
14 & 13,363 & 21,024 & \frac{13,363}{114,688}\approx 0.11652 & \frac{21,024}{114,688}\approx 0.18331 \\
&&&&\\
15 & 28,635 & 49,856 & \frac{28,635}{245,760}\approx 0.11652 & \frac{49,856}{245,760}\approx 0.20286 \\
&&&&\\
16 & 61,088 & 108,976 & \frac{61,088}{524,288}\approx 0.11652 & \frac{108,976}{524,288}\approx 0.20786 \\
&&&&\\
17 & 129,812 & 224,976 & \frac{129,812}{1,114,112}\approx 0.11652 & \frac{224,976}{1,114,112}\approx 0.20193 \\
&&&&\\
18 & 274,896 & 517,552 & \frac{274,896}{2,359,296}\approx 0.11652 & \frac{517,552}{2,359,296}\approx 0.21937 \\
&&&&\\
19 & 580,336 & 1,111,856 & \frac{580,336}{4,980,736}\approx 0.11652 & \frac{1,111,856}{4,980,736}\approx 0.22323 \\
&&&&\\
20 & 1,221,760 & 2,273,680 & \frac{1,221,760}{10,485,760}\approx 0.11652 & \frac{2,273,680}{10,485,760}\approx 0.21684 \\
&&&&\\
21 & 2,565,696 & 5,124,736 & \frac{2,565,696}{22,020,096}\approx 0.11652 & \frac{5,124,736}{22,020,096}\approx 0.23273 \\
&&&&\\
22 & 5,375,744 & 10,879,712 & \frac{5,375,744}{46,137,344}\approx 0.11652 & \frac{10,879,712}{46,137,344}\approx 0.23581 \\
&&&&\\
23 & 11,240,192 & 22,105,536 & \frac{11,240,192}{96,468,992}\approx 0.11652 & \frac{22,105,536}{96,468,992}\approx 0.22915 \\
&&&&\\
24 & 23,457,792 & 49,096,752 & \frac{23,457,792}{201,326,592}\approx 0.11652 & \frac{49,096,752}{201,326,592}\approx 0.24387 \\
&&&&\\
25 & 48,870,400 & 103,338,816 & \frac{48,870,400}{419,430,400}\approx 0.11652 & \frac{103,338,816}{419,430,400}\approx 0.24638 \\
&&&&\\
26 & 101,650,432 & 208,999,264 & \frac{101,650,432}{872,415,232}\approx 0.11652 & \frac{208,999,264}{872,415,232}\approx 0.23956 \\
&&&&\\
27 & 211,120,128 & 459,059,616 & \frac{211,120,128}{1,811,939,328}\approx 0.11652 & \frac{459,059,616}{1,811,939,328}\approx 0.25335 \\
&&&&\\
\hline
\end{array}$$

We can also take the construction of a $Q_2$-free subgraph of $Q_5$ and divide the non-edges into $e$ and $o$ edges so that each $Q_3$ contains at least one $e$-edge and one $o$-edge.  This can be done as follows:

$$\begin{array}{|c|c|c|c|}
  \hline
   \text{$e$ edges} & \begin{array}{cc}
                \text{[*0001]} & \text{[*1000]} \\
                \text{[*0111]} & \text{[*1110]} \\
                \text{[1010*]} & \text{[1101*]} \\
                \text{[0*101]} & \text{[0*010]} \\
                \text{[10*10]} & \text{[11*01]} \\
                \text{[001*1]} & \text{[010*1]}
              \end{array}
   &
  \text{$o$ edges} &  \begin{array}{cc}
                \text{[*0100]} & \text{[*0010]} \\
                \text{[*1101]} & \text{[*1011]} \\
                \text{[0000*]} & \text{[0111*]} \\
                \text{[1*000]} & \text{[1*111]} \\
                \text{[00*11]} & \text{[01*00]} \\
                \text{[111*0]} & \text{[100*1]}
              \end{array} \\
  \hline
\end{array}$$

Since there are $80-24=56$ $a$ edges and $24$ non-$a$ edges.  This yields the recursion: $$e_{k+4} = 16e_k + 40p_k + 24\cdot 2^{k-2}.$$  This will in general give worse bounds than the $aeo$-colored $Q_4$, but it gives a better result for $k=5$ since $p_5=0$.  This gives $c(Q_3,9) \le 320$ which will in turn give smaller values for $c(Q_3,9+3k)$.

We may also be able to improve some other small bounds by taking a $Q_2$-free $Q_m$ and then dividing up the non-edges into $e$- and $o$-edges.  For instance, it would give $c(Q_3,10)\le 736$ which seems to be the last number where improvements would be possible in this manner.

\section{A General Construction for $Q_3$-free subgraphs of the hypercube}
We denote an edge as before in the form $[x_1x_2\cdots x_{i-1}*x_{i+1}\cdots x_n]$.  Where $x_j\in \{0,1\}$.  We denote two function for an edge:

$p(e)$ is the number of ones before the $*$ minus the number of ones afterward.

Then let $$A=\{e: p(e)\equiv 0\pmod 4\}.$$

For the second part, we could consider instead edges where $p(e)\equiv 1,2,\text{ or } 3\pmod 4$ and the argument would still apply.  Hence $A$ contains at most 1/4 of the edges of $Q_n$.

It remains to show that $A$ contains at least one edge from each $Q_3$ in $Q_n$.  Denote a $Q_3$ by $a*b*c*d$ where $a,b,c,d$ are strings of zeros and ones.  Let $|s|$ denote the number of ones in a string $s$.

\emph{Case 1:} $|a|+|b|+|c|+|d|\equiv 0 \pmod 2$.
Then consider the edges:$$a*b0c0d\qquad \text{and} \qquad a*b1c1d.$$

Both of them are part of our $Q_3$ and one of them has $p(e)\equiv 0$.

\emph{Case 2:} $|a|+|b|+|c|+|d|\equiv 1 \pmod 2$.
Then consider the edges:$$a1b*c0d\qquad \text{and} \qquad a0b*c1d.$$

Again, both of them are part of our $Q_3$ and one of them has $p(e)\equiv 0$.  Hence $Q_n- A$ contains no $Q_3$.

\end{document}